\documentclass{amsart}

\theoremstyle{definition}

\theoremstyle{remark}

\reversemarginpar

\begin{document}

\title[Beta function]{The integrals in Gradshteyn and Ryzhik. \\ 
Part 6: The beta function.}

\author{Victor H. Moll}
\address{Department of Mathematics,
Tulane University, New Orleans, LA 70118}
\email{vhm@math.tulane.edu}

\subjclass{Primary 33}

\date{\today}


\begin{abstract}
We present a systematic derivation of some of the 
definite integrals in the classical 
table of Gradshteyn and Ryzhik that can be reduced to the beta function.
\end{abstract}

\maketitle

\newcommand{\nn}{\nonumber}
\newcommand{\ba}{\begin{eqnarray}}
\newcommand{\ea}{\end{eqnarray}}
\newcommand{\ift}{\int_{0}^{\infty}}
\newcommand{\ione}{\int_{0}^{1}}
\newcommand{\ifft}{\int_{- \infty}^{\infty}}
\newcommand{\no}{\noindent}
\newcommand{\realpart}{\mathop{\rm Re}\nolimits}
\newcommand{\imagpart}{\mathop{\rm Im}\nolimits}

\newtheorem{Definition}{\bf Definition}[section]
\newtheorem{Thm}[Definition]{\bf Theorem} 
\newtheorem{Example}[Definition]{\bf Example} 
\newtheorem{Lem}[Definition]{\bf Lemma} 
\newtheorem{Note}[Definition]{\bf Note} 
\newtheorem{Cor}[Definition]{\bf Corollary} 
\newtheorem{Prop}[Definition]{\bf Proposition} 
\newtheorem{Problem}[Definition]{\bf Problem} 
\numberwithin{equation}{section}

\section{Introduction} \label{intro} 
\setcounter{equation}{0}

The table of integrals \cite{gr} contains some evaluations that can be 
derived by elementary means from the {\em beta function}, defined by 
\begin{equation}
B(a,b) = \int_{0}^{1} x^{a-1}(1-x)^{b-1} \, dx.
\label{beta-def}
\end{equation}
\noindent
The convergence of the integral in (\ref{beta-def}) requires $a, \, b > 0$. 
This definition appears as $\mathbf{3.191.3}$ in \cite{gr}.  

Our goal is to present in a systematic manner, the evaluations appearing in 
the classical table of Gradshteyn and Ryzhik \cite{gr}, that involve this 
function. In this part, we restrict to algebraic integrands leaving the 
trigonometric forms for a future publication. This paper 
complements \cite{moll-gr4} that dealt with the 
{\em gamma function} defined by
\begin{equation}
\Gamma(a) := \ift x^{a-1}e^{-x} \, dx. 
\end{equation}
\noindent
These functions are related by the functional equation
\begin{equation}
B(a,b) = \frac{\Gamma(a) \, \Gamma(b)}{\Gamma(a+b)}. 
\label{functional}
\end{equation}
\noindent
A proof of this identity can be found in \cite{irrbook}.  \\

The special values $\Gamma(n) = (n-1)!$ and 
\begin{equation}
\Gamma \left( n + \tfrac{1}{2} \right) = \frac{\sqrt{\pi}}{2^{2n}} 
\, \frac{(2n)!}{n!} 
\label{gamma1a}
\end{equation}
\noindent
for $n \in \mathbb{N}$, will 
be used to simplify the values of the integrals presented here. Proofs of these 
formulas can be found in \cite{moll-gr4} as  well as in Proposition 
\ref{value12} below. 

The other property that will be employed frequently is 
\begin{equation}
\Gamma(a) \, \Gamma(1-a) = \frac{\pi}{\sin \pi a}.
\label{reflection}
\end{equation}
\noindent
The reader will find in \cite{irrbook} a proof based on the product 
representation of these functions. A challenging problem is to produce a
proof that only employs changes of variables. 

The table \cite{gr} contains some direct values: 
\begin{equation}
\int_{0}^{1} \frac{x^{p} \, dx}{(1-x)^{p}} = \frac{p \pi}{\sin p \pi}
\end{equation}
\noindent
is $\mathbf{3.192.1}$ and is 
evaluated by identifying it as $B(p+1,1-p)$. Formula 
$\mathbf{3.192.2}$ is 
\begin{equation}
\int_{0}^{1} \frac{x^{p} \, dx}{(1-x)^{p+1}} = - \frac{\pi}{\sin p \pi} 
\end{equation}
\noindent
has the value $B(p+1,-p) = \Gamma(p+1)\Gamma(-p)$. Next, $\mathbf{3.192.3}$ is 
\begin{equation}
\int_{0}^{1} \frac{(1-x)^{p}}{x^{p+1}} \, dx = - \frac{\pi}{\sin p \pi}
\end{equation}
\noindent
and the change of variables $t = 1/x$ in $\mathbf{3.192.4}$  produces
\begin{equation}
\int_{1}^{\infty} (x-1)^{p-1/2} \frac{dx}{x} = 
\int_{0}^{1} t^{-p-1/2} (1-t)^{p-1/2} \, dt
\end{equation}
\noindent
and this is 
\begin{equation}
B \left( \tfrac{1}{2}-p, \tfrac{1}{2} +p \right) = 
\Gamma \left( \tfrac{1}{2} -p \right) 
\Gamma \left( \tfrac{1}{2} +p \right)  = \frac{\pi}{\cos p \pi}, 
\end{equation}
\noindent
as stated in \cite{gr}.

Let $b = \tfrac{1}{2}$ in (\ref{beta-def}) to obtain
\begin{equation}
\int_{0}^{1} \frac{x^{a-1} \, dx}{\sqrt{1-x}} = B \left(a, \tfrac{1}{2} \right) 
= \frac{\Gamma(a) \, \sqrt{\pi}}{\Gamma \left( a + \tfrac{1}{2} \right)}. 
\end{equation}
\noindent
The special values $a = n+1$ and $a=n + \tfrac{1}{2}$ appear as 
$\mathbf{3.226.1}$
and $\mathbf{3.226.2}$, respectively.

\section{Elementary properties} \label{sec-elemen} 
\setcounter{equation}{0}

Many of the properties of the beta function can be established by simple 
changes of variables. For example, letting $y = 1-x$ in (\ref{beta-def}) 
yields the symmetry 
\begin{equation}
B(a,b) = B(b,a).
\end{equation}

It should  not be surprising that a clever change of variables might lead 
to a beautiful result. This is illustrated following Serret \cite{serret1}. 
Start with 
\begin{eqnarray}
B(a,a) & = & \int_{0}^{1} (x - x^{2})^{a-1} \, dx  \nonumber \\
 & = & 2 \int_{0}^{1/2} \left[ \tfrac{1}{4} 
- \left( \tfrac{1}{2} - x \right)^{2} \right]^{a-1} \, dx.  \nonumber 
\end{eqnarray}
\noindent
The natural change of variables $v = \tfrac{1}{2}-x$ yields 
\begin{equation}
B(a,a) = 2 \int_{0}^{1/2} \left( \tfrac{1}{4} - v^{2} \right)^{a-1} \, dv. 
\end{equation}
\noindent
The next step is now clear: let  $s = 4v^{2}$ to produce 
\begin{equation}
B(a,a) = 2^{1-2a} B \left(a, \tfrac{1}{2} \right). 
\end{equation}
\noindent
The functional equation (\ref{functional}) converts this identity into 
Legendre's original form:

\begin{Prop}
\label{value12}
The gamma function satisfies 
\begin{equation}
\Gamma \left( a + \tfrac{1}{2}  \right) = \frac{\Gamma(2a) \,
\Gamma( \tfrac{1}{2} ) }{\Gamma(a) \, 2^{2a-1}}. 
\label{Legendre}
\end{equation}
\noindent
In particular, for $a=n \in \mathbb{N}$, this yields (\ref{gamma1a}).
\end{Prop}

\section{Elementary changes of variables} \label{sec-elchanges} 
\setcounter{equation}{0}

The integral (\ref{beta-def}) defining the beta function can be transformed by
changes of variables. For example, the new variable $x = t/u$, reduces 
(\ref{beta-def}) to 
\begin{equation}
\int_{0}^{u} t^{a-1}(u-t)^{b-1} \, dt = u^{a+b-1}B(a,b),
\end{equation}
\noindent
that appears as 
$\mathbf{3.191.1}$ in \cite{gr}. The effect of this change of variables
is to express the beta function as an integral over a finite interval. Observe
that the integrand vanishes at both end points.  Similarly, the change 
$t = (v-u)x+u$ maps the interval $[0,1]$ to $[u,v]$. It yields 
\begin{equation}
\int_{u}^{v} (t-u)^{a-1} (v-t)^{b-1} \, dt = (v-u)^{a+b-1} B(a,b). 
\end{equation}
\noindent
This is $\mathbf{3.196.3}$ in \cite{gr}. The special case $u=0, \, v=n$ and 
$a = \nu, \, b = n+1$ appears as $\mathbf{3.193}$ in \cite{gr} as 
\begin{equation}
\int_{0}^{n} x^{\nu-1} (n-x)^{n} \, dx = \frac{n^{\nu+n} \, n!}
{\nu(\nu+1)(\nu+2) \cdots (\nu+n)}. 
\end{equation}

\medskip

Several integrals in \cite{gr} can be obtained by a small variation of the 
definition. For example, the integral
\begin{equation}
\int_{0}^{1} ( 1 - x^{a})^{b-1} \, dx = \frac{1}{a} B \left(1/a, 
b \right)
\end{equation}
\noindent
can be obtained by the change of variables $t = x^{a}$. This appears as 
$\mathbf{3.249.7}$ in 
\cite{gr} and illustrates the fact that it not necessary for 
the integrand to vanish at {\em both} end points. The special case 
$a=2$ appears as $\mathbf{3.249.5}$:
\begin{equation}
\int_{0}^{1} (1-x^{2})^{b-1} \, dx  = \tfrac{1}{2} B \left( \tfrac{1}{2}, 
b \right) = 2^{2b-2}B(b,b), 
\end{equation}
\noindent
where the second identity follows from Legendre's duplication formula
(\ref{Legendre}). \\

The change of variables 
$t = cx$ produces a scaled version:
\begin{equation}
\int_{0}^{c} (c^{a}-t^{a})^{b-1} \, dt = \frac{1}{a}c^{a(b-1)+1}
B \left(1/a, b \right). 
\end{equation}
\noindent
The special case $a=2$ yields 
\begin{equation}
\int_{0}^{c} (c^{2}-t^{2})^{b-1} \, dt = \frac{c^{2b-1}}{2}
B \left( 1/2,  b \right). 
\end{equation}
\noindent
The choice $b = n + \tfrac{1}{2}$ appears as 
$\mathbf{3.249.2}$ in \cite{gr}:
\begin{equation}
\int_{0}^{c} (c^{2}-t^{2})^{n - 1/2} \, dt = \frac{\pi c^{2n}}{2^{2n+1}} 
\binom{2n}{n}. 
\end{equation}

\medskip 

Similarly
$\mathbf{3.251.1}$ in \cite{gr} is 
\begin{equation}
\int_{0}^{1} x^{c-1}(1-x^{a})^{b-1} \, dx = \frac{1}{a} B \left( \frac{c}{a},
b \right).
\end{equation}

\medskip

The change of variables $t = 1/x$ converts (\ref{beta-def}) into
\begin{equation}
\int_{1}^{\infty} t^{-a-b} (t-1)^{b-1} \, dt = B(a,b).
\end{equation}
\noindent
Letting $t = x^{p}$ yields 
\begin{equation}
\int_{1}^{\infty} x^{p(1-a-b)-1} \left(x^{p}-1 \right)^{b-1} \, dx = 
\frac{1}{p}B(a,b).
\end{equation}
\noindent
The special case $\nu = b $ and $\mu = p(1-a-b)$ is 
$\mathbf{3.251.3}$:
\begin{equation}
\int_{1}^{\infty} x^{\mu-1} \left(x^{p}-1 \right)^{\nu-1} \, dx = 
\frac{1}{p} B \left( 1 - \nu - \mu/p, \nu \right).
\end{equation}

\medskip

\section{Integrals over a half-line} \label{sec-halfline} 
\setcounter{equation}{0}

The beta function can also be expressed as an integral over a half-line. The 
change of variables $t = x/(1-x)$ maps $[0,1]$ onto $[0, \infty)$ and it 
produces from (\ref{beta-def})
\begin{equation}
B(a,b) = \ift \frac{t^{a-1} \, dt}{(1+t)^{a+b}}. 
\label{beta-halfline}
\end{equation}
\noindent
In particular, if $a+b=1$, using (\ref{functional}) and (\ref{reflection}), 
we obtain
\begin{equation}
\ift \frac{t^{a-1} \, dt}{1+t} = \frac{\pi}{\sin \pi a}. 
\label{euler-90}
\end{equation}
\noindent
This can be scaled to produce, for $a>0$ and $c>0$,
\begin{equation}
\ift \frac{x^{a-1} \, dx}{x+c} = \frac{\pi}{\sin \pi a } c^{a-1}
\quad \text{ for } c>0
\label{int55}
\end{equation}
\noindent
that appears as $\mathbf{3.222.2}$ in \cite{gr}.  In 
the case $c<0$ we have a singular
integral. Define $b = -c >0$ and $s = x/b$, so now we have to evaluate
\begin{equation}
I = -b^{a-1} \ift \frac{s^{a-1} \, ds}{1-s}.
\end{equation}
\noindent
The integral is considered as a Cauchy principal value
\begin{equation}
I = \lim\limits_{\epsilon \to 0} \int_{0}^{1} \frac{s^{a-1} \, ds}
{(1-s)^{1-\epsilon}} + \int_{1}^{\infty} 
\frac{s^{a-1} \, ds}{(1-s)^{1-\epsilon}}. 
\end{equation}
\noindent
Let $y = 1/s$ in the second integral and evaluate them in terms of the 
beta function to produce
\begin{equation}
I = \lim\limits_{\epsilon \to 0} 
\epsilon \Gamma(\epsilon) \times 
\frac{1}{\epsilon} \left( \frac{\Gamma(a)}{\Gamma(a+ \epsilon)} -
\frac{\Gamma(1-a-\epsilon)}{ \Gamma(1-a)} \right).
\end{equation}
\noindent
Use L'Hopital's rule to evaluate and obtain
\begin{equation}
I = -\frac{\Gamma'(a)}{\Gamma(a)} + \frac{\Gamma'(1-a)}{\Gamma(a)}.
\end{equation}
\noindent
Using the relation $\Gamma(a) \Gamma(1-a) = \pi \text{cosec }\pi a$, this 
reduces to $\pi \cot \pi a$. Therefore we have
\begin{equation}
\ift \frac{x^{a-1} \, dx}{x+c} = -\frac{\pi}{\tan \pi a } (-c)^{a-1}
\quad \text{ for } c<0
\label{int56}
\end{equation}
\noindent
The change of variables $x = e^{-t}$ produces, for $c<0$,
\begin{equation}
\int_{-\infty}^{\infty} \frac{e^{-\mu t} \, dt}{e^{-t}+c} = 
-\pi \cot(\mu \pi) \, (-c)^{\mu-1}. 
\end{equation}
\noindent
The special case $c=-1$ appears as $\mathbf{3.313.1}$:
\begin{equation}
\int_{-\infty}^{\infty} \frac{e^{-\mu t} \, dt}{1- e^{-t}} = 
\pi \cot(\mu \pi). 
\end{equation}

\medskip

We now consider several examples in \cite{gr} that are direct 
consequences of (\ref{int55}) and (\ref{int56}). In the first example, we
combine (\ref{int55}) with the partial fraction decomposition
\begin{equation}
\frac{1}{(x+a)(x+b)} = \frac{1}{b-a} \left( \frac{1}{x+a} - \frac{1}{x+b}
\right)
\end{equation}
\noindent
leads to $\mathbf{3.223.1}$:
\begin{equation}
\ift \frac{x^{\mu-1} \, dx}{(x+b)(x+a)} = 
\frac{\pi}{b-a} ( a^{\mu-1} - b^{\mu-1}) \text{cosec}(\pi \mu).
\end{equation}
\noindent
Similarly, 
\begin{equation}
\frac{1}{x+b} - \frac{1}{x-a} = \frac{a+b}{(a-x)(b+x)}
\end{equation}
\noindent
leads to $\mathbf{3.223.2}$:
\begin{equation}
\ift \frac{x^{\mu-1} \, dx}{(b+x)(a-x)} = \frac{\pi}{a+b} 
\left( b^{\mu-1} \text{ cosec}(\mu \pi) + 
a^{\mu-1} \cot(\mu \pi) \right),
\end{equation}
\noindent
using (\ref{int55}) and (\ref{int56}). The result $\mathbf{3.223.3}$:
\begin{equation}
\ift \frac{x^{\mu-1} \, dx}{(a-x)(b-x)} = \pi \cot(\mu \pi) 
\frac{a^{\mu-1}-b^{\mu-1}}{b-a},
\end{equation}
\noindent
follows from 
\begin{equation}
\frac{1}{(a-x)(b-x)} = \frac{1}{a-b} \left( \frac{1}{b-x}-\frac{1}{a-x} 
\right). 
\end{equation}
\noindent
Finally, $\mathbf{3.224}$:
\begin{equation}
\ift \frac{(x+b)x^{\mu-1} \, dx}{(x+a)(x+c)} = \frac{\pi}{\sin(\mu \pi)}
\left( \frac{a-b}{a-c} a^{\mu-1} + \frac{c-b}{c-a} c^{\mu-1} \right),
\end{equation}
\noindent
follows from 
\begin{equation}
\frac{x+b}{(x+a)(x+c)} = \frac{b-a}{c-a} \frac{1}{x+a} -
                         \frac{b-c}{c-a} \frac{1}{x+c}. 
\end{equation}

\medskip

We can now transform (\ref{beta-halfline}) to the interval $[0,1]$ by splitting 
$[0,\infty)$ as $[0,1]$ followed by $[1,\infty)$. In the second integral, we 
let $t = 1/s$. The final result is 
\begin{equation}
B(a,b) = \int_{0}^{1} \frac{t^{a-1} + t^{b-1}}{(1+t)^{a+b}} \, dt.
\label{beta-12}
\end{equation}
\noindent
This formula, that appears as $\mathbf{3.216.1}$, makes 
it apparent that the beta
function is symmetric: $B(a,b) = B(b,a)$. The change of variables $s = 1/t$
converts (\ref{beta-12}) into $\mathbf{3.216.2}$:
\begin{equation}
B(a,b) = \int_{1}^{\infty} \frac{s^{a-1} + s^{b-1}}{(1+s)^{a+b}} \, ds.
\label{beta-13}
\end{equation}
\noindent
It is easy to introduce a parameter: let $c>0$ and consider the 
change of variables $t = cx$ in (\ref{beta-halfline}) 
to obtain
\begin{equation}
\ift \frac{x^{a-1} \, dx}{(1+cx)^{a+b}} = c^{-a} B(a,b), 
\label{beta-halfline1}
\end{equation}
\noindent
that appears as $\mathbf{3.194.3}$. We 
can now shift the lower limit of integration via
$t = x+u$ to produce 
\begin{equation}
\int_{u}^{\infty} (t-u)^{a-1} (t + v)^{-a-b} \, dt = 
(u + v)^{-b} B(a,b), 
\end{equation}
\noindent
where $v = 1/c-u$. This is $\mathbf{3.196.2}$, where $v$ is denoted by
$\beta$. Now let $b = c-a$ in the 
special case $v = 0$ to obtain
\begin{equation}
\int_{u}^{\infty} (t-u)^{a-1}t^{-c} \, dt = u^{a-c}B(a,c-a). 
\end{equation}
\noindent
This appears as $\mathbf{3.191.2}$. \\

We now write (\ref{beta-halfline}) using the change of  variables 
$t = x^{c}$. It produces 
\begin{equation}
\ift \frac{x^{ac-1} \, dx}{(1+x^{c})^{a+b}} = \frac{1}{c}B(a,b).
\end{equation}
\noindent
The special case $c=2$ and $a = 1 + \mu/2, \, b = 1 - \mu/2$ produces
$\mathbf{3.251.6}$ in the form
\begin{equation}
\ift \frac{x^{\mu+1} \, dx}{(1+x^{2})^{2}} = \frac{\mu \pi}
{4 \, \sin \mu \pi/2}.
\end{equation}

Now let $b = 1-a$ and choose $a = p/c$ to obtain
\begin{equation}
\ift \frac{x^{p-1} \, dx}{1+x^{c}} = \frac{1}{c} B \left( \frac{p}{c}, 
\frac{c-p}{c} \right) = \frac{\pi}{c} \text{ cosec}(\pi p/c). 
\end{equation}
\noindent
This appears as $\mathbf{3.241.2}$ in \cite{gr}. \\

Similar arguments establish $\mathbf{3.196.4}$:
\begin{equation}
\int_{1}^{\infty} \frac{dx}{(a-bx)(x-1)^{\nu}} = - \frac{\pi}{b} 
\text{ cosec}(\nu \pi) \left( \frac{b}{b-a} \right)^{\nu}.
\end{equation}
\noindent
Indeed, the change of variables $t = x-1$ yields
\begin{equation}
\int_{1}^{\infty} \frac{dx}{(a-bx)(x-1)^{\nu}} = 
\ift \frac{dt}{\left[ (a-b) -bt \right] \, t^{\nu}}, 
\end{equation}
\noindent
and scaling via the new variable $z = bt/(b-a)$ gives
\begin{equation}
\int_{1}^{\infty} \frac{dx}{(a-bx)(x-1)^{\nu}} = 
- \frac{1}{b} \left( \frac{b}{b-a} \right)^{\nu} 
\ift \frac{dz}{(1+z) \, z^{\nu}}. 
\end{equation}
\noindent
The result follows from (\ref{beta-halfline}) and the value
\begin{equation}
B(\nu, 1- \nu) = \Gamma(\nu) \Gamma(1- \nu) = \frac{\pi}{\sin \pi \nu}. 
\end{equation}
\noindent
The same argument gives $\mathbf{3.196.5}$:

\begin{equation}
\int_{-\infty}^{1} \frac{dx}{(a-bx)(1-x)^{\nu}} =  \frac{\pi}{b} 
\text{ cosec}(\nu \pi) \left( \frac{b}{a-b} \right)^{\nu}.
\end{equation}
\noindent

\section{Some direct evaluations} \label{sec-direct}
\setcounter{equation}{0}

There are many more integrals in \cite{gr} that can be evaluated in terms
of the beta function. For example, $\mathbf{3.221.1}$ states that
\begin{equation}
\int_{a}^{\infty} \frac{(x-a)^{p-1} \, dx}{x-b} = \pi (a-b)^{p-1} 
\text{ cosec }\pi p. 
\end{equation}
\noindent
To establish these identities, we assume 
that $a>b$ to avoid the singularities. The change of variables 
$t = (x-a)/(a-b)$ yields
\begin{equation}
\int_{a}^{\infty} \frac{(x-a)^{p-1} \, dx}{x-b} = 
(a-b)^{p-1} \ift \frac{t^{p-1} \, dt}{1+t},
\end{equation}
\noindent
and this integral appears in (\ref{euler-90}).

Similarly, $\mathbf{3.221.2}$ states that
\begin{equation}
\int_{-\infty}^{a} \frac{(a-x)^{p-1} \, dx}{x-b} = 
-\pi (b-a)^{p-1} \text{ cosec }\pi p. 
\end{equation}
\noindent
This is evaluated by the change of variables $y=(a-x)/(b-a)$.

\medskip

The table contains several evaluations that are elementary corollaries of 
(\ref{beta-halfline}). Starting with 
\begin{equation}
\ift \frac{x^{a} \, dx }{(1+x)^{b}} = B(a+1,b-a-1) = 
\frac{\Gamma(a+1) \, \Gamma(b-a-1)}{\Gamma(b)},
\label{new-beta}
\end{equation}
\noindent
we find the case $a=p$ and $b=3$ in $\mathbf{3.225.3}$:
\begin{equation}
\ift \frac{x^{p} \, dx}{(1+x)^{3}} = \frac{\Gamma(p+1) \, \Gamma(2-p)}
{\Gamma(3)} = \frac{p(1-p)}{2} \frac{\pi}{\sin(p \pi)},
\end{equation}
\noindent
using elementary properties of the gamma function.  \\

The change of variables $t=1+x$ converts (\ref{new-beta}) into 
\begin{equation}
\int_{1}^{\infty} \frac{(t-1)^{a} \, dt }{t^{b}} = B(a+1,b-a-1) = 
\frac{\Gamma(a+1) \, \Gamma(b-a-1)}{\Gamma(b)}.
\label{new-beta1}
\end{equation}
\noindent
The special case $a=p-1$ and $b=2$ gives
\begin{equation}
\int_{1}^{\infty} \frac{(t-1)^{p-1} \, dt}{t^{2}} = \Gamma(p) \Gamma(2-p) 
= (1-p) \Gamma(p) \Gamma(1-p) = \frac{\pi (1-p)}{\sin(p \pi)}. 
\end{equation}
\noindent
This appears as $\mathbf{3.225.1}$. Similarly, the 
case $a=1-p$ and $b=3$ produces
$\mathbf{3.225.2}$:

\begin{equation}
\int_{1}^{\infty} \frac{(t-1)^{1-p} \, dt}{t^{3}} = 
\frac{\Gamma(2-p) \Gamma(1+p) }{\Gamma(3)} 
= \frac{1}{2}p(1-p) \Gamma(p) \Gamma(1-p) = 
\frac{\pi \, p (1-p)}{2\sin(p \pi)}. 
\end{equation}

\section{Introducing parameters} \label{sec-param} 
\setcounter{equation}{0}

It is often convenient to introduce free parameters in a definite integral.
Starting with (\ref{beta-halfline}), the change of variables 
$t = \frac{u}{v}x^{c}$ yields
\begin{equation}
B(a,b) = cu^{a}v^{b} \ift \frac{t^{ac-1} \, dt}{(v + ut^{c})^{a+b}}. 
\label{master}
\end{equation}
\noindent
This formula appears as $\mathbf{3.241.4}$ in \cite{gr} with the parameters 
\begin{equation}
a = \frac{\mu}{\nu}, \, b = n+1 - \frac{\mu}{\nu}, \, c = \nu, \, u = q, 
\text{ and } v = p, 
\end{equation}
\noindent
in the form
\begin{equation}
\ift \frac{x^{\mu-1} \, dx}{(p + q x^{\nu})^{n+1}} = 
\frac{1}{\nu \, p^{n+1}} \left( \frac{p}{q} \right)^{\mu/\nu} \, 
\frac{\Gamma(\mu/\nu) \, \Gamma(n+1- \mu/\nu)}{\Gamma(n+1)}. 
\nonumber
\end{equation}
\noindent
This is a messy notation and it leaves the wrong impression that $n$ should 
be an integer. \\

\noindent
$\bullet$ The special case $v=c=1$ and $b=p+1-a$ produces
\begin{equation}
\ift \frac{t^{a-1} \, dt}{(1+ut)^{p+1}} = \frac{1}{u^{a}} B(a,p+1-a).
\label{beta-11}
\end{equation}
\noindent
This appears as $\mathbf{3.194.4}$ 
in \cite{gr}, except that it is written in terms
of binomial coefficients as
\begin{equation}
\ift \frac{t^{a-1} \, dt}{(1+ut)^{p+1}} = 
(-1)^{p} \frac{\pi}{u^{a}} \binom{a-1}{p} \, \text{cosec}(\pi a).
\end{equation}
\noindent
We prefer the notation in (\ref{beta-11}). \\

\noindent
$\bullet$ The special case $v=c=1$ and $b=2-a$ produces
\begin{equation}
\ift \frac{t^{a-1} \, dt}{(1+ut)^{2}} = \frac{1}{u^{a}} B(a,2-a). 
\end{equation}
\noindent
Using (\ref{functional}) and (\ref{reflection}) yields the form
\begin{equation}
\ift \frac{t^{a-1} \, dt}{(1+ut)^{2}} = 
\frac{(1-a) \, \pi}{u^{a} \, \sin \pi a}. 
\end{equation}
\noindent
This appears as $\mathbf{3.194.6}$ in \cite{gr}. \\

\noindent
$\bullet$ The special case $u=v=1$ and $c=q$, and choosing $a = p/q$ and 
$b=2-p/q$ yields $\mathbf{3.241.5}$ in the form
\begin{equation}
\ift \frac{x^{p-1} \, dx}{(1+x^{q})^{2}} = \frac{q-p}{q^{2}} 
\, \frac{\pi}{\sin(\pi p/q)}.
\end{equation}

\medskip

\noindent
$\bullet$ The special case $c=1$ and $a = m+1, \, b = n -m - \tfrac{1}{2}$
produces 
\begin{equation}
\ift \frac{t^{m} \, dt}{(v + ut)^{n + \tfrac{1}{2} }} = 
\frac{1}{u^{m+1} \, v^{n-m-\tfrac{1}{2}}}  
B \left( m+1, n-m- \tfrac{1}{2} \right)
\end{equation}
\noindent
Using (\ref{functional}) and (\ref{gamma1a}) this reduces to
\begin{equation}
\ift \frac{t^{m} \, dt}{(v + ut)^{n + \tfrac{1}{2} }} = 
\frac{m! \, n! \, (2n-2m-2)!}{(n-m-1)! \, (2n)!} 2^{2m+2}
\frac{v^{m-n+1/2}}{u^{m+1}}, 
\end{equation}
\noindent
for $m, \, n \in \mathbb{N}$, with $n>m$. This 
appears as $\mathbf{3.194.7}$ in \cite{gr}. \\

\noindent
$\bullet$ The special case $u=v=1$ and $ b = \tfrac{1}{2}-a$ yields
\begin{equation}
\ift \frac{t^{ac-1} \, dt}{\sqrt{1+t^{c} } } = 
\frac{1}{c} B \left(a, \tfrac{1}{2}-a \right). 
\end{equation}
\noindent
Writing $a = p/c$ we recover $\mathbf{3.248.1}$: 
\begin{equation}
\ift \frac{t^{p-1} \, dt}{\sqrt{1+t^{c} } } = 
\frac{1}{c} B \left( \tfrac{p}{c}, \tfrac{1}{2}- \tfrac{p}{c} \right). 
\end{equation}

\noindent

$\bullet$ Now replace $v$ by $v^{2}$ in (\ref{master}). Then, with $u=1, \, 
a = \tfrac{1}{2}, \, c = 2$, so that $ac=1$ and $b=n- \tfrac{1}{2}$ we 
obtain
\begin{equation}
\ift \frac{dt}{(v^{2}+t^{2})^{n}} = \frac{1}{2v^{2n-1}} B \left( 
\tfrac{1}{2}, n - \tfrac{1}{2} \right). 
\end{equation}
\noindent
This can be written as 
\begin{equation}
\ift \frac{dt}{(v^{2}+t^{2})^{n}} = 
\frac{\sqrt{\pi} \, \Gamma(n - 1/2)}{2 \Gamma(n) v^{2n-1}}
\end{equation}
\noindent
that appears as $\mathbf{3.249.1}$ in \cite{gr}. \\

$\bullet$ The special case $v=1, \, c=2$ and $b = \frac{n}{2}-a$ in 
(\ref{master}) yields
\begin{equation}
\ift \frac{t^{2a-1} \, dt}{(1 + ut^{2})^{n/2}} = 
\frac{1}{2u^{a}} B \left( a, \tfrac{n}{2} - a \right). 
\end{equation}
\noindent
Now $a = 1/2$ gives
\begin{equation}
\ift ( 1 + ut^{2})^{-n/2} \, dt  =  \frac{1}{2 \sqrt{u}} 
B \left( \tfrac{1}{2}, \tfrac{n-1}{2} \right)
 =  \frac{\sqrt{\pi}}{2 \sqrt{u}} 
\frac{\Gamma(\tfrac{n-1}{2})}{\Gamma(n/2)}. 
\end{equation}
\noindent
It is curious that the table \cite{gr} contains $\mathbf{3.249.8}$ as
the special case $u = 1/(n-1)$ of this evaluation.  \\

$\bullet$ We now put  $u=v=1$ and $c=2$ in (\ref{master}). Then, with
$b = 1-\nu -a$ and $a = \mu/2$, we obtain $\mathbf{3.251.2}$: 
\begin{equation}
\ift \frac{t^{\mu -1} \, dt}{(1+t^{2})^{1- \nu}} = 
\frac{1}{2} B \left( \frac{\mu}{2}, 1 - \nu - \frac{\mu}{2} \right).
\end{equation}

\medskip

$\bullet$ We now consider the case $c=2$ in (\ref{master}): 
\begin{equation}
\ift \frac{t^{2a-1} \, dt}{(v + ut^{2})^{a+b}} = \frac{1}{2u^{a}v^{b}} 
B(a,b). 
\end{equation}
\noindent
The special case $a = m + \tfrac{1}{2}$ and $b = n-m + \tfrac{1}{2}$ yields
\begin{equation}
\ift \frac{t^{2m} \, dt}{(v+ut^{2})^{n+1}} = 
\frac{\Gamma(m+1/2) \, \Gamma(n-m+1/2)}{2u^{m+1/2} v^{n-m+1/2} \Gamma(n+1)}, 
\end{equation}
\noindent
and using (\ref{gamma1a}) we obtain $\mathbf{3.251.4}$:
\begin{equation}
\ift \frac{t^{2m} \, dt}{(v+ut^{2})^{n+1}} = 
\frac{\pi (2m)! (2n-2m)!}{2^{2n+1} m! (n-m)! n! \, u^{m + 1/2} v^{n-m+1/2}}, 
\end{equation}
\noindent
for $n, m \in \mathbb{N}$ with $n > m$. 

On the other hand, if we choose $a = m+1$ and $b = n-m$ we obtain 
$\mathbf{3.251.5}$:
\begin{equation}
\ift \frac{t^{2m+1} \, dt}{(v+ut^{2})^{n+1}} = 
\frac{\Gamma(m+1) \, \Gamma(n-m)}{2u^{m+1} v^{n-m} \Gamma(n+1)} = 
\frac{m! (n-m-1)!}{2n! u^{m+1}v^{n-m}}. 
\end{equation}

\medskip

Several evaluation in \cite{gr} come from the form
\begin{equation}
\int_{0}^{1} t^{aq-1} (1-t^{q})^{b-1} \, dt = \frac{1}{q} B(a,b), 
\label{621}
\end{equation}
\noindent
obtained from (\ref{beta-def}) by the change of variables $x = t^{q}$. \\

$\bullet$ The choice $a = 1 + p/q$ and $b = 1 -p/q$ produces 
\begin{equation}
\int_{0}^{1} t^{p+q-1}(1-t^{q})^{-p/q} \, dt = 
\frac{1}{q} B \left( 1 + \frac{p}{q}, 1-\frac{p}{q} \right) =
\frac{p \pi}{q^{2}} \text{ cosec} \left( \frac{p \pi}{q} \right). 
\end{equation}
\noindent
This appears as $\mathbf{3.251.8}$. \\

$\bullet$ The choice $a = 1/p$ and $b = 1 -1/p$ gives
\begin{equation}
\int_{0}^{1} x^{q/p-1} ( 1 - x^{q})^{-1/p} \, dx = 
\frac{1}{q} B \left( \frac{1}{p}, 1 - \frac{1}{p} \right) = 
\frac{\pi}{q} \text{ cosec} \left( \frac{ \pi}{p} \right). 
\end{equation}
\noindent
This appears as $\mathbf{3.251.9}$. \\

$\bullet$ The reader can now check that the choice $a = p/q$ and $b=1-p/q$ 
yields the evaluation
\begin{equation}
\int_{0}^{1} x^{p-1} ( 1 - x^{q})^{-p/q} \, dx = 
\frac{1}{q} B \left( \frac{p}{q}, 1 - \frac{p}{q} \right) = 
\frac{\pi}{q} \text{ cosec} \left( \frac{ p \pi}{q} \right). 
\end{equation}
\noindent
This appears as $\mathbf{3.251.10}$. \\

$\bullet$ Putting $v=1$ and $b = \nu -a$ in (\ref{master}) we get 
\begin{equation}
\ift \frac{t^{ac-1} \, dt}{(1+ut^{c})^{\nu}} = \frac{1}{cu^{a}} B(a, \nu -a).
\end{equation}
\noindent
Now let $a = r/c$ to obtain
\begin{equation}
\ift \frac{t^{r-1} \, dt}{(1+ut^{c})^{\nu}} = \frac{1}{cu^{r/c}} 
B \left( \frac{r}{c}, \nu - \frac{r}{c} \right).
\end{equation}
\noindent
This appears as $\mathbf{3.251.11}$. \\

$\bullet$ We now choose $b=1-1/q$ in (\ref{621}) to obtain
\begin{equation}
\int_{0}^{1} \frac{t^{aq-1} \, dt}{ \sqrt[q]{1-t^{q}}} = 
\frac{1}{q}B \left(a,1 - \frac{1}{q} \right).
\end{equation}
\noindent
Finally, writing $a = c - (m-1)/q$ gives the form
\begin{equation}
\int_{0}^{1} \frac{t^{cq-m} \, dt}{ \sqrt[q]{1-t^{q}}} = 
\frac{1}{q}B \left(c + \frac{1}{q}-\frac{m}{q},1 - \frac{1}{q} \right).
\end{equation}
\noindent
The special case $q=2$ produces
\begin{equation}
\int_{0}^{1} \frac{t^{2c-m} \, dt}{\sqrt{1-t^{2}}} = 
\tfrac{1}{2}B \left( c + \tfrac{1}{2}-\tfrac{m}{2},  \tfrac{1}{2} \right) 
= \frac{\Gamma(c + \tfrac{1}{2} - \tfrac{m}{2}) \sqrt{\pi}}
{2 \Gamma(c+1 - \tfrac{m}{2})}. 
\end{equation}
\noindent
In particular, if $c = n+1$ and $m=1$ we obtain $\mathbf{3.248.2}$:
\begin{equation}
\int_{0}^{1} \frac{t^{2n+1} \, dt}{\sqrt{1-t^{2}}} = 
\frac{\sqrt{\pi} \, n!}{2 \Gamma(n + 3/2)} = 
\frac{2^{2n} \, n!^{2}}{(2n+1)!}. 
\end{equation}
\noindent
Similarly, $c=n$ and $m=0$ yield $\mathbf{3.248.3}$:
\begin{equation}
\int_{0}^{1} \frac{t^{2n} \, dt}{\sqrt{1-t^{2}}} = \frac{\pi}{2^{2n+1}}  \, 
\frac{(2n)!}{n!^{2}} = \frac{\pi}{2^{2n+1}} \binom{2n}{n}. 
\end{equation}

\medskip

In the case $q=3$ we get
\begin{equation}
\int_{0}^{1} \frac{t^{3c-m} \, dt}{ \sqrt[3]{1-t^{3}}} = 
\frac{1}{3}B \left(c + \frac{1}{3}-\frac{m}{3},1 - \frac{1}{3} \right).
\end{equation}
\noindent
This includes $\mathbf{3.267.1}$ and $\mathbf{3.267.2}$ in \cite{gr}:
\begin{eqnarray}
\int_{0}^{1} \frac{t^{3n} \, dt}{\sqrt[3]{1-t^{3}}} & =  & 
\frac{2 \pi}{3 \sqrt{3}} \frac{\Gamma(n + \tfrac{1}{3})}{\Gamma(\tfrac{1}{3}) 
\, \Gamma( n+1)} \nonumber \\
\int_{0}^{1} \frac{t^{3n-1} \, dt}{\sqrt[3]{1-t^{3}}} & =  & 
 \frac{(n-1)! \Gamma( \tfrac{2}{3})}{3 \Gamma(n + \tfrac{2}{3}) } 
\nonumber
\end{eqnarray}
\noindent
The latest edition of \cite{gr} has added our suggestion
\begin{equation}
\int_{0}^{1} \frac{t^{3n-2} \, dt}{\sqrt[3]{1-t^{3}}}  =  
 \frac{\Gamma(n- \tfrac{1}{3}) \, \Gamma( \tfrac{2}{3})}
{3 \Gamma(n + \tfrac{1}{3}) } 
\end{equation}
\noindent
as $\mathbf{3.267.3}$. 

\section{The exponential scale} \label{sec-expon} 
\setcounter{equation}{0}

We now present examples of (\ref{beta-def}) written in terms of the 
exponential function. The change of variables $x = e^{-ct}$ in (\ref{beta-def})
yields
\begin{equation}
\ift e^{-at} ( 1 - e^{-ct})^{b-1} \, dt = \frac{1}{c} B \left( \frac{a}{c},b 
\right). 
\end{equation}
\noindent
This appears as $\mathbf{3.312.1}$ in \cite{gr}. On the other hand, if we let 
$x = e^{-ct}$ in (\ref{beta-halfline})  we get 
\begin{equation}
\int_{-\infty}^{\infty} \frac{e^{-act} \, dt}{(1+e^{-ct})^{a+b}} = 
\frac{1}{c} B(a,b).
\end{equation}
\noindent
This appears as $\mathbf{3.313.2}$ in \cite{gr}. The reader can now use the 
techniques described above to verify
\begin{equation}
\int_{-\infty}^{\infty} \frac{e^{- \mu x} \, dx }{(e^{b/a} + 
e^{-x/a})^{\nu}} = a \, 
\text{exp} \, \left[ b \left( \mu - \frac{\nu}{a} \right)
\right] \, B \left( a \mu, \nu - a \mu \right),
\end{equation}
\noindent
that appears as $\mathbf{3.314}$. The 
choice $b=0, \, \nu=1$ and relabelling parameters
by $a = 1/q$ and $\mu=p$ yields $\mathbf{3.311.3}$:
\begin{equation}
\int_{-\infty}^{\infty} \frac{e^{-px} \, dx}{1+ e^{-qx}} = 
\frac{1}{q} B \left( \frac{p}{q}, 1 - \frac{p}{q} \right) = 
\frac{\pi}{q} \text{cosec}\left(\frac{\pi p}{q} \right),
\end{equation}
\noindent
using the identity $B(x,1-x) = \pi \text{cosec}(\pi x)$ in the last step.
This is the form given in the table.  \\

The integral $\mathbf{3.311.9}$:
\begin{equation}
\int_{- \infty}^{\infty} 
 \frac{e^{-\mu x} \, dx}{b+e^{-x}} = \pi b^{\mu-1} \text{ cosec}(\mu \pi)
\end{equation}
\noindent
can be evaluated via the change of variables $t = e^{-x}/b$ and 
(\ref{euler-90}) to produce
\begin{equation}
I = b^{\mu-1} \ift \frac{t^{\mu-1} \, dt}{1+t}. 
\end{equation}

\section{Some logarithmic examples} \label{sec-logar} 
\setcounter{equation}{0}

The beta function appears in the evaluation of definite integrals involving
logarithms. For example, $\mathbf{4.273}$ states that
\begin{equation}
\int_{u}^{v} \left( \ln \frac{x}{u} \right)^{p-1} \, 
\left( \ln \frac{v}{x} \right)^{q-1} \, \frac{dx}{x} = 
B(p,q) \, \left( \ln \frac{v}{u} \right)^{p+q-1}. 
\end{equation}
\noindent
The evaluation is simple: the change of variables $x = ut$ produces, with 
$c = v/u$, 
\begin{equation}
I = \int_{1}^{c} \ln^{p-1}t \, ( \ln c - \ln t)^{q-1} \, \frac{dt}{t}.
\end{equation}
\noindent
The change of variables $z = \frac{\ln t}{\ln c}$ give the result.  \\

A second example is $\mathbf{4.275.1}$:
\begin{equation}
\int_{0}^{1} \left[ ( - \ln x)^{q-1} - x^{p-1}(1-x)^{q-1} \right] \, dx 
= \frac{\Gamma(q)}{\Gamma(p+q)} \left[ \Gamma(p+q) - \Gamma(p) \right], 
\end{equation}
\noindent
that should be written as 
\begin{equation}
\int_{0}^{1} \left[ ( - \ln x)^{q-1} - x^{p-1}(1-x)^{q-1} \right] \, dx 
= \Gamma(q) -  B(p,q). 
\end{equation}
\noindent
The evaluation is elementary, using Euler form of the gamma function
\begin{equation}
\Gamma(q) = \int_{0}^{1} \left( - \ln x \right)^{q-1} \, dx.
\end{equation}

\section{Examples with a fake parameter} \label{sec-fake} 
\setcounter{equation}{0}

The evaluation $\mathbf{3.217}$:
\begin{equation}
\ift \left( \frac{b^{p}x^{p-1}}{(1+bx)^{p}} - 
\frac{(1+bx)^{p-1}}{b^{p-1}x^{p}} \right) \, dx = \pi \cot \pi p
\end{equation}
\noindent
has the obvious parameter $b$. We say that this is a {\em fake parameter} in 
the sense that a simple scaling shows that the integral is independent of it.
Indeed, the change of variables 
$t = bx$ shows this independence. Therefore the evaluation amounts to
showing  that
\begin{equation}
\ift \left( \frac{t^{p-1}}{(1+t)^{p}} - 
\frac{(1+t)^{p-1}}{t^{p}} \right) \, dt = \pi \cot \pi p.
\label{nine-two}
\end{equation}
\noindent
To achieve this, we let $y = 1/t$ in the second integral to produce
\begin{equation}
\lim\limits_{\epsilon \to 0} 
\ift \frac{t^{p-1- \epsilon} \, dt }{(1+t)^{p}} -
\ift \frac{t^{\epsilon-1} \, dt }{(1+t)^{1-p}}. 
\end{equation}
\noindent
The integrals above evaluate to $B(p-\epsilon,\epsilon)-B(\epsilon,1-p-
\epsilon)$. Using 
\begin{equation}
B(a,b) = \frac{\Gamma(a) \Gamma(b)}{\Gamma(a+b)} \text{ and }
\Gamma(a) \Gamma(1-a) = \frac{\pi}{\sin(\pi a)}
\end{equation}
\noindent
this reduces to 
\begin{equation}
I = \lim\limits_{\epsilon \to 0} 
\epsilon \Gamma(\epsilon) 
\left( 
\frac{
\Gamma(p-\epsilon) \Gamma(p+\epsilon) \sin(\pi(p+ \epsilon)) - 
\Gamma^{2}(p) \sin(\pi p)}
{\epsilon \Gamma(p) \Gamma(p+\epsilon) \sin(\pi(p+\epsilon))}
\right).
\end{equation}
\noindent
Now recall that 
\begin{equation}
\lim\limits_{\epsilon \to 0} \epsilon \Gamma(\epsilon) =1
\end{equation}
\noindent
and reduce the previous limit to 
\begin{equation}
I = \frac{1}{\Gamma^{2}(p) \, \sin(\pi p)} 
\lim\limits_{\epsilon \to 0} 
\frac{1}{\epsilon} \left( \Gamma(p-\epsilon) \Gamma(p+\epsilon) 
\sin(\pi(p+\epsilon)) -
\Gamma^{2}(p) \sin(\pi p) \right). 
\end{equation}
Using L'Hopital's rule we find that $I = \pi \cot(\pi p)$ as required. \\

The example $\mathbf{3.218}$ 
\begin{equation}
\ift \frac{x^{2p-1}-(a+x)^{2p-1}}{(a+x)^{p} \, x^{p}} \, dx = 
\pi \cot \pi p
\end{equation}
\noindent
also shows a fake parameter. The change of variable $x = at$ reduces the 
integral above to
\begin{equation}
\ift \frac{t^{2p-1}-(1+t)^{2p-1}}{(1+t)^{p} \, t^{p}} \, dt = 
\pi \cot \pi p
\end{equation}
\noindent
This can be written as 
\begin{equation}
I = \ift 
\left(\frac{t^{p-1}}{(1+t)^{p}} - \frac{(1+t)^{p-1}}{t^{p}} \right) \, dt.
\end{equation}
\noindent
The result now follows from (\ref{nine-two}). \\

\section{Another type of logarithmic integral} \label{sec-newtype} 
\setcounter{equation}{0}

Entry $\mathbf{4.251.1}$ is
\begin{equation}
\ift \frac{x^{a-1} \, \ln x}{x+b} \, dx = \frac{\pi \, b^{a-1}}{\sin \pi a} 
\left( \ln b - \pi \cot \pi a \right). 
\end{equation}

\noindent
To check this evaluation we first scale by $x = bt$ and obtain
\begin{equation}
\ift \frac{x^{a-1} \, \ln x}{x+b} \, dx = 
b^{a-1} \ln b \ift \frac{t^{a-1} \,dt}{1+t} + b^{a-1} \ift \frac{t^{a-1} \, 
\ln t}{1+t} \, dt.
\end{equation}
\noindent
The first integral is simply 
\begin{equation}
\ift \frac{t^{a-1} \, dt }{1+t} = B(a,1-a) = \Gamma(a) \Gamma(1-a) = 
\frac{\pi}{\sin \pi a}.
\end{equation}
\noindent
The second one is evaluated as 
\begin{equation}
\ift \frac{t^{a-1} \, \ln t}{1+t} \, dt = - \pi^{2} \frac{\cos \pi a}
{\sin^{2}(\pi a )} 
\end{equation}
\noindent
by differentiating (\ref{beta-halfline}) 
with respect to $a$. The evaluation follows from here.

\section{A hyperbolic looking integral} \label{sec-hyperbolic} 
\setcounter{equation}{0}

The evaluation of $\mathbf{3.457.3}$:
\begin{equation}
\int_{-\infty}^{\infty} \frac{x \, dx}{(a^{2}e^{x} + e^{-x})^{\mu}} = 
-\frac{1}{2a^{\mu}} B \left( \frac{\mu}{2}, \frac{\mu}{2} \right) \, \ln a,
\end{equation}
\noindent
is done as follows: write
\begin{equation}
I = \frac{1}{a^{\mu}} 
\int_{-\infty}^{\infty} \frac{x \, dx}{(ae^{x} + a^{-1}e^{-x})^{\mu}} 
\end{equation}
\noindent
and let $t = ae^{x}$ to produce
\begin{equation}
I = \frac{1}{a^{\mu}} \ift \frac{t^{\mu-1} \, (\ln t - \ln a) \, dt}
{(1+t^{2})^{\mu}}.
\end{equation}
\noindent
The change of variables $s = t^{2}$ yields
\begin{equation}
I = \frac{1}{4a^{\mu}} \ift \frac{s^{\mu/2-1} \, \ln s \, ds}{(1+s)^{\mu}} -
\frac{\ln a}{2a^{\mu}} \ift \frac{s^{\mu/2-1} \, ds}{(1+s)^{\mu}}. 
\end{equation}
\noindent
The first integral vanishes. This follows directly from the change $s \mapsto
1/s$. The second integral is the beta value indicated in the formula. 

In particular, the value $a=1$ yields
\begin{equation}
\int_{-\infty}^{\infty} \frac{x \, dx}{\cosh^{\mu}x} = 0.
\end{equation}
\noindent
Differentiating with respect to $\mu$ produces
\begin{equation}
\int_{-\infty}^{\infty}  x \ln \cosh x \, dx = 0,
\end{equation}
\noindent
that appears as $\mathbf{4.321.1}$ in \cite{gr}.

\bigskip

\noindent
{\bf Acknowledgments}. The author wishes to thank Luis Medina for a 
careful reading of the manuscript. The partial support of 
\noindent
$\text{NSF-DMS } 0409968$ is also acknowledged.

\end{document}